\documentclass[11pt]{amsart}
\usepackage{amssymb}
\usepackage{amsxtra}
\theoremstyle{plain}

\newcommand{\cleqn}{\setcounter{equation}{0}}
\newcommand{\clth}{\setcounter{theorem}{0}}
\newcommand {\sectionnew}[1]{\section{#1}\cleqn\clth}

\newtheorem{theorem}{Theorem}[section]
\newtheorem{lemma}[theorem]{Lemma}
\newtheorem{definition-theorem}[theorem]{Definition-Theorem}
\newtheorem{proposition}[theorem]{Proposition}
\newtheorem{corollary}[theorem]{Corollary}
\newtheorem{definition}[theorem]{Definition}
\newtheorem{example}[theorem]{Example}
\newtheorem{remark}[theorem]{Remark}
\newtheorem{notation}[theorem]{Notation}
\newcommand \bth[1] { \begin{theorem}\label{t#1} }
\newcommand \ble[1] { \begin{lemma}\label{l#1} }

\newcommand \bpr[1] { \begin{proposition}\label{p#1} }
\newcommand \bco[1] { \begin{corollary}\label{c#1} }
\newcommand \bde[1] { \begin{definition}\label{d#1}\rm }
\newcommand \bex[1] { \begin{example}\label{e#1}\rm }
\newcommand \bre[1] { \begin{remark}\label{r#1}\rm }

\newcommand \bnota[1] { \begin{notation}\label{n#1}\rm }
\renewcommand {\eth} { \end{theorem} }
\newcommand {\ele} { \end{lemma} }

\newcommand {\epr} { \end{proposition} }
\newcommand {\eco} { \end{corollary} }
\newcommand {\ede} { \end{definition} }
\newcommand {\eex} { \end{example} }
\newcommand {\ere} { \end{remark} }

\newcommand {\enota} { \end{notation} }
\newcommand \thref[1]{Theorem \ref{t#1}}
\newcommand \leref[1]{Lemma \ref{l#1}}
\newcommand \prref[1]{Proposition \ref{p#1}}
\newcommand \coref[1]{Corollary \ref{c#1}}
\newcommand \deref[1]{Definition \ref{d#1}}

\newcommand \lb[1]{\label{#1}}
\def \d {{\partial}}   

\def \Rset {{\mathbb R}}         
\def \Cset {{\mathbb C}}

\def \DD  {{\mathcal{D}}}

\def \LL  {{\mathcal{L}}}
\def \De {\Delta}   

\def \al {\alpha}

\def \Ga {\Gamma}

\def \mt  {\mapsto}

\def \hra {\hookrightarrow}

\def \ci  {\circ}

\def \rcor {\rangle}
\def \lcor {\langle}

\def \ol {\overline}
\def \wt {\widetilde}

\def \id { {\mathrm{id}} }

\def \Ad { {\mathrm{Ad}} }

\def \Lie { {\mathrm{Lie \,}} }
\def \g  {\mathfrak{g}}   

\def \h  {\mathfrak{h}}
\def \d  {\mathfrak{d}}

\def \n  {\mathfrak{n}}

\def \b  {\mathfrak{b}}

\def \q  {\mathfrak{q}}

\def \uu {\mathfrak{u}}
\def \tt {\mathfrak{t}}
\def \l {\mathfrak{l}}

\renewcommand \Im { {\mathrm{Im}} }

\begin{document}
\title[Quotients of Drinfeld and Heisenberg doubles]
{Weak splittings of quotients of \\ Drinfeld and Heisenberg doubles}
\author[Milen Yakimov]{Milen Yakimov}
\address{Department of Mathematics \\
Louisiana State Univerity \\
Baton Rouge, LA 70803 \\
U.S.A.
}
\email{yakimov@math.lsu.edu}
\thanks{The author was supported in part
by NSF grants DMS-1001632 and and DMS-1303036.}
\date{}
\begin{abstract} 
We investigate the fine structure of the 
simplectic foliations of Poisson homogeneous spaces. 
Two general results are proved for weak splittings of surjective 
Poisson submersions from Heisenberg and Drinfeld doubles. 
The implications of these results are that the torus 
orbits of symplectic leaves of the quotients 
can be explicitly realized as Poisson--Dirac submanifolds
of the torus orbits of the doubles. The results have a wide 
range of applications to many families of real and complex 
Poisson structures on flag varieties.
Their torus orbits of leaves recover important families of 
varieties such as the open Richardson varieties. 
\end{abstract}
\maketitle
\sectionnew{Introduction}\lb{Introduction}
The geometry of Poisson--Lie groups is well understood, 
both in the case of the standard Poisson structures on complex simple Lie groups 
\cite{HL,KZ} and the general Belavin--Drinfeld \cite{BD} Poisson structures 
\cite{Y1}. The torus orbits of symplectic leaves in the former case 
are the double Bruhat cells of the simple Lie group. One of the fundamental results in 
the theory of cluster algebras is the 
Berenstein--Fomin--Zelevinsky theorem \cite{BFZ}
that their coordinate rings are upper cluster algebras.
Recently, the coordinate rings of the 
${\mathrm{SL}}_n$ groups, equipped with the Cremmer--Gervais Poisson structures 
from \cite{BD}, were also shown to be upper cluster algebras \cite{GSV}. 
The motivation for these results is that cluster algebras give rise to 
Poisson structures by the work of Gekhtman, Shapiro, and Vainshtein \cite{GSV0}, 
and one attempts to go in the opposite direction using Poisson varieties from 
the theory of quantum groups.

On the other hand, the possible cluster algebra structures on coordinate rings 
of torus orbits of symplectic leaves of Poisson homogeneous spaces is much less well
understood. In the special case of the standard complex Poisson structures on 
flag varieties, this is precisely the problem of constructing 
cluster algebra structures on the coordinate rings of the open Richardson varieties.
These varieties have been recently studied in \cite{BC,KLS,RW}
in relation to Schubert calculus and total positivity. Chevalier \cite{Ch} conjectured 
cluster algebra structures for the Richardson strata in the case when one of the two Weyl group elements
is a parabolic Coxeter element. Leclerc \cite{L} generalized this construction
and showed that the coordinate ring of each open Richardson variety contains 
a cluster algebra whose rank is equal to the dimension of the variety. 
These cluster structures come from an additive categorification.
Another cluster algebra structure on the Richardson--Lusztig 
strata of the Grassmannian was conjectured by Muller and Speyer \cite{MS}
using Postnikov diagrams \cite{Po}. 

In this paper we prove a very general result that realizes torus orbits 
of symplectic leaves of a large class of Poisson homogeneous 
spaces as Poisson--Dirac submanifolds of torus orbits of symplectic leaves of 
Drinfeld and Heisenberg doubles. 
It applies to many important families of complex and real Poisson structures 
on flag varieties, double flag varieties, and their generalizations.
In the context of cluster algebras, the point of this construction is that 
the coordinate rings of affine Poisson varieties with conjectured cluster algebra 
structures are realized as quotients of better understood coordinate rings of 
Poisson varieties,
some of which are already proven to possess cluster algebra structures.
The ideals defining these quotients are not Poisson but have somewhat similar properties 
coming from a notion of ``weak splitting of surjective Poisson submersions.''
The construction of the latter is the main point of the paper.   

To explain this in precise terms, we recall that 
to each point of a Poisson homogeneous space of a Poisson--Lie group
Drinfeld \cite{Dr} associated a Lagrangian subalgebra of the double 
and proved an equivariance property of this map. Motivated 
by this construction, Lu and Evens associated to each 
quadratic Lie algebra (a complex or real Lie algebra 
$\d$ equipped with a nondegenerate invariant symmetric bilinear 
form $\lcor.,.\rcor$) the variety of its Lagrangian subalgebras 
$\LL(\d, \lcor.,.\rcor)$ and initiated its systematic study in 
\cite{EL1}. This is a singular projective variety. 

Fix a connected Lie group $D$ with Lie algebra $\d$. Given any pair 
of Lagrangian subalgebras $\g_\pm$ such that $\d$ is the vector space direct sum 
of $\g_+$ and $\g_-$ (in other words, given a Manin triple $(\d, \g_+, \g_-)$
with respect to the bilinear form $\lcor.,.\rcor$), one defines the $r$-matrix 
$r = \frac{1}{2}\sum_j \xi_j \wedge x_j$ where $\{x_j\}$ and $\{\xi_j\}$ is a
pair of dual bases of $\g_+$ and $\g_-$. Using the adjoint action of $D$ on 
$\LL(\d, \lcor.,.\rcor)$, we construct 
the bivector field $\chi(r)$ on $\LL(\d, \lcor.,.\rcor)$.
Here and below $\chi$ refers to the infinitesimal action associated
to a Lie group action. It was proved in \cite{EL1} that $\chi(r)$ is Poisson.
Up to minor technical details, the singular projective Poisson variety 
\[
(\LL(\d, \lcor.,.\rcor), \chi(r))
\]
captures the geometry of all Poisson homogeneous spaces of the 
Poisson--Lie groups integrating the Lie bialgebras $\g_\pm$. 
The $D$-orbits on $\LL(\d, \lcor.,.\rcor)$ are 
compete Poisson submanifolds, i.e., they are unions of symplectic leaves of $\chi(r)$. 
They have the form  $D/N(\l)$ where $\l$ is a Lagrangian 
subalgebra of $(\d, \lcor.,.\rcor)$ and $N(\l)$ is the normalizer of $\l$ in $D$.
The rank of the Poisson structure $\chi(r)$ at each point of $D/N(\l)$ is given by  
\cite[Theorem 4.10]{LY}. This describes the coarse structure of 
the symplectic foliations of the spaces $(D/N(\l), \chi(r))$ or equivalently 
the variety of Lagrangian subalgebras $(\LL(\d, \lcor.,.\rcor), \chi(r))$.

In this paper we address the problem of describing the fine structure of the symplectic 
foliations of these spaces. From the point of view of 
Lie theory and cluster algebras, the most important $D$-orbits in this picture
are the orbits
\[
(D/N(\g_+), \chi(r)) \hra ( \LL(\d, \lcor.,.\rcor), \chi(r)).
\]  
These Poisson varieties capture all examples of real and complex 
Poisson structures on flag varieties and double flag varieties that 
appeared in previous studies, see e.g. \cite{EL1,FL,GY,WY}.
The Poisson varieties $(D/N(\g_+), \chi(r))$ also have the properties 
that they are quotients of {\em{Drinfeld and Heisenberg doubles}}. 
Recall that those are the Poisson varieties
\[
(D, \pi= L(r) - R(r)) \quad \mbox{and} \quad 
(D, \pi'=L(r) + R(r)),
\]
respectively. Here and below $R(.)$ and $L(.)$ refer to right and left 
invariant multivector fields on 
a Lie group. The Poisson structure $\pi$ vanishes along the group
\[
H:= N(\g_+) \cap N(\g_-)
\]
and as a consequence the left and right action of 
$H$ on $D$ preserves both Poisson structures $\pi$ and $\pi'$. The corresponding 
Poisson reductions will be denoted by 
\[
(D/H, \pi_H) \; \; \mbox{and} \; \; (D/H, \pi'_H).
\]
The canonical projections
\begin{multline}
\mu \colon (D/H, \pi_H) \to (D/N(\g_+), \chi(r)) 
\; \; \mbox{and} \\ 
\mu' \colon (D/H, \pi'_H) \to (D/N(\g_+), \chi(r))
\label{subm}
\end{multline}
are {\em{surjective Poisson submersions}}.

We prove that under certain general assumptions the symplectic 
leaves of $(D/N(\g_+), \chi(r))$ can be realized as 
explicit symplectic submanifolds of the symplectic leaves
of the (reduced) {\em{Drinfeld double}} $(D/H, \pi_H)$ or {\em{Heisenberg double}} 
$(D/H, \pi'_H)$ (or even both in some cases). 
To be more precise, recall \cite{CF}
that a Poisson manifold $(X, \pi_X)$ is {\em{a Poisson--Dirac submanifold
admitting a Dirac projection}} of a Poisson manifold $(M, \Pi)$ 
if $X$ is a submanifold of $M$, and there exists a subbundle $E$ of 
$T_XM$ such that 
\[
E \oplus TX = T_X M \; \; \mbox{and} \; \;
\Pi- \pi_X \in \Ga(X, \wedge^2 E).
\]
In this framework we find an explicit construction 
of sections of the surjective Poisson submersion 
$\mu \colon (D/H, \pi_H) \to (D/N(\g_+), \chi(r))$
over each $N(\g_-)$-orbit whose images are 
Poisson--Dirac submanifolds of $(D/H, \pi_H)$.
This is the {\em{weak splitting 
of the first surjective Poisson 
submersion}} in \eqref{subm} from a Drinfeld double.
(We refer the reader to Sect. \ref{sect2} below 
and \cite[Sect. 2]{GY} for the definition of the notion 
and additional background.)
As a corollary of the general construction, 
the symplectic leaves within each $N(\g_-)$-orbit
on $(D/N(\g_+), \pi_{D/N(\g_+)})$ are uniformly 
embedded as symplectic submanifolds of $(D/H, \pi_H)$.
Similarly, we  
construct sections of the second surjective Poisson submersion
$\mu' \colon (D/H, \pi'_H) \to (D/N(\g_+), \chi(r))$ in \eqref{subm}
over each $N(\g_+)$-orbit whose images are 
Poisson--Dirac submanifolds of $(D/H, \pi'_H)$ admitting a Dirac projection.
These constructions of weak splittings work under 
certain general assumptions, see Theorems \ref{tmain2} 
and \ref{tmain2H} for details. In Sect. \ref{sect5}
we show that the conditions are satisfied for many important 
families of Poisson structures.

The above results have a wide range of applications. 
In the case of the standard Poisson structures on 
flag varieties we recover the weak splittings from \cite{GY}. 
Double flag varieties arise naturally as closed strata in 
partitions of wonderful group compactifications \cite{DP}. This 
gives rise to Poisson structures on them that are not 
products of Poisson structures on each factor
\cite{WY}. The second splitting result above for Heisenberg doubles 
is applicable to this family of Poisson varieties. The real forms of a complex simple
Lie algebra $\g$ give rise to real Poisson structures 
on the related complex flag variety defined in \cite{FL}. 
Again the above second splitting is applicable for this family. 
Finally, the Delorme's classification result in \cite{De} gives rise 
to canonical Poisson structures on products of flag varieties for complex 
simple Lie groups (i.e, flag varieties for a reductive group).
Except for some very special cases, these 
Poisson structures are not products of Poisson structures on the 
factors. Our weak splittings are applicable for those families too. 

Another motivation for the results in the paper is the study of the 
spectra of the quantizations of the homogeneous coordinate rings 
of the above mentioned families of varieties. Currently,
only the spectra of quantum flag varieties are understood \cite{Y2}. 
We expect that a quantum version of our weak splittings of surjective 
Poisson submersions will be helpful in understanding the spectra of the 
quantizations of these families of varieties on the basis of the 
works on spectra of quantum groups \cite{HL,J}. It appears that 
such quantum weak splittings should be also closely related to the 
notion of quantum folding of Berenstein and Greenstein \cite{BG}.

The paper is organized as follows. In Sect. \ref{sect2} we 
review the notion of Poisson--Dirac submanifolds of Poisson 
manifolds, and weak sections and weak splittings of surjective 
Poisson submersions. In Sect. \ref{sect3} we prove two general 
theorems on the construction of weak sections and weak splittings 
for quotients of Drinfeld doubles. In Sect. \ref{sect4} 
similar theorems are proved for quotients of Heisenberg doubles.
Sect. \ref{sect5} contains applications of these theorems.

We finish the introduction with some notation that will be used throughout 
the paper. Given a group $G$, $d \in G$ and two subgroups $H_1$ and
$H_2$ of $G$, we will denote the $H_1$-orbit through $d H_2 \in G/H_2$
by 
\[
H_1 \cdot d H_2 \subset G/H_2
\]
(to distinguish it from the double coset $H_1 d H_2 \subset G$).

For a Lie group $G$, we will denote by $G^\circ$ 
its identity component. For a smooth manifold $X$, $X^\circ$ will 
denote a connected component of $X$. Given a Lie group $G$ 
and a subalgebra $\uu$ of its Lie algebra $\g$, we will 
denote by $N(\uu)$ the normalized of $\uu$ in $G$ with respect
to the adjoint action. Finally, recall that a Poisson structure $\pi$ on a manifold 
$M$ gives rise to the bundle map $\pi^\sharp \colon T^* M \to T M$, given by 
\[
\pi^\sharp(\al) = \al \otimes \id(\pi), \; \; \al \in T_m M, m \in M.
\] 
\\ \hfill \\
{\bf Acknowledgements.} I would like to thank Bernard Leclerc for his valuable 
suggestions on the first version of the preprint. I am also grateful to 
I. Dimitrov, G. Mason, S. Montgomery, I. Penkov, V. S. Varadarajan, and 
J. Wolf for the opportunity to present these and related results 
at Lie theory meetings at UCB, UCLA, UCSC, and USC.
\sectionnew{Poisson--Dirac submanifolds and weak splittings}
\label{sect2}
\noindent
In this section we review the notion of weak splitting of a 
surjective Poisson submersion from \cite{GY}.
We start by recalling several facts about Poisson--Dirac submanifolds
of Poisson manifolds, introduced and studied in \cite{CF,X, V}.

\bde{PD} A submanifold $X$ of a Poisson manifold 
$(M, \Pi)$ is called a {\em{Poisson--Dirac submanifold}} if the following 
conditions are satisfied:

(i) For each symplectic leaf $S$ of $(M, \Pi)$,
the intersection $S\cap X$ is clean {\rm(}i.e., it is smooth and
$T_x (S \cap X)= T_x S \cap T_x X$
for all $x \in S \cap X${\rm)}, and
$S\cap X$ is a symplectic submanifold of
$(S, (\Pi|_S)^{-1})$.

(ii) The family of symplectic structures
$(\Pi|_S)^{-1}|_{S\cap X}$ is induced by a smooth Poisson
structure $\pi$ on $X$.
\ede

Clearly, in the setting of \deref{PD}, the symplectic leaves 
of $(X, \pi)$ are the connected components of the intersections
of symplectic leaves of $(M, \pi)$ with $X$.

An important criterion is provided by the following 
result.

\bpr{PD2} {\rm{[Crainic, Fernandes]}} \cite{CF} Assume that 
$X$ is a submanifold of a Poisson manifold $(M, \Pi)$ for which 
there exits a subbundle $E \subset T_X M$ such that 

(i)  $E \oplus T X = T_X M$ and

(ii) $\Pi \in \Ga(X, \wedge^2 T X \oplus \wedge^2 E).$
\\
Then $X$ is Poisson--Dirac submanifold of $(M, \Pi)$.
\epr

In the setting of \prref{PD2} the projection of 
$\Pi|_X$ into $\Ga(X, \wedge^2 TX)$ along $\wedge^2 E$ 
is exactly the needed 
Poisson structure $\pi$ in \deref{PD}.  
Poisson--Dirac submanifolds with the property of \prref{PD2}
are called {\em{Poisson--Dirac submanifolds admitting
a Dirac projection}} by Crainic--Fernandes \cite{CF} and quasi-Poisson
submanifolds by Vaisman \cite{V}. We will use the former 
term.

In the presence of the condition (i) in \prref{PD2}, the condition (ii) 
is equivalent \cite{CF} to 
\begin{equation}
\label{iii}
\Pi_m^\sharp ((T_m X)^0) \subset E_m, \quad \forall m \in X.
\end{equation}
Here and below for a subspace $V \subseteq T_m M$, $V^0$ will denote its 
orthogonal complement in $T_m^* M$.

We continue with the notions of weak sections and weak splittings of surjective Poisson submersions.

\bde{wsdef} Assume that $(M, \Pi)$ and $(N, \pi)$ are Poisson
manifolds, $X$ is a Poisson submanifold of $(N, \pi)$, 
and that $p \colon (M, \Pi) \to (N, \pi)$ is a surjective
Poisson submersion. A {\em{weak section}} of $p$ over 
$X$ is a smooth map $i \colon X \to M$ such that
$p \circ i = \id_X$ and $i(X)$ is a Poisson--Dirac
submanifold of $(M, \Pi)$ with induced Poisson 
structure $i_*(\pi|_X)$.
\ede

In this situation we derive from \prref{PD2} an
explicit realization of all symplectic leaves of $(X, \pi|_X)$
in terms of those of $(M, \Pi)$:
  
{\em{In the setting of \deref{wsdef} one has that each symplectic
leaf of $(X, \pi|_X)$ has the form 
$i^{-1}((X \cap S)^\ci)$ where $S$ is a symplectic leaf of
$(M, \Pi)$. In addition $i$ realizes explicitly all leaves 
of $(X, \pi)$ as symplectic submanifolds of symplectic leaves 
of $(M, \Pi)$.}} 

\bre{alg} The following special case of the notion of weak section 
has an equivalent algebraic characterization 
which is of particular interest, see \cite[Proposition 2.6]{GY} for 
details.

Let $p \colon (M, \Pi) \to (N, \pi)$ be a surjective
Poisson submersion, $X$ be an open subset of $N$, and  
$i \colon X \to M$ be a smooth map such that
$p \circ i = \id_X$ and $i(X)$ is a smooth submanifold 
of $N$. In particular, $p^* : (C^\infty(N), \{.,.\}_\pi) \to (C^\infty(M), \{.,.\}_\Pi)$ 
is a homomorphism of Poisson algebras. Then $i$ is a weak section 
with associated bundle equal to the tangent bundle to the fibers of $p$ if and only if 
\[
i^* : (C^\infty(M), \{.,.\}_\Pi) \to (C^\infty(N), \{.,.\}_\pi)
\]
is a homomorphism of Poisson $(C^\infty(N), \{.,.\}_\pi)$-modules 
with respect to the action on the first term coming from $p^*$.
\ere

\bde{weakspldef} Let $(M, \Pi)$ and $(N, \pi)$ be Poisson
manifolds and $p \colon (M, \Pi) \to (N, \pi)$ be a surjective
Poisson submersion. A {\em{weak splitting}} of $p$ is a partition
\[
N = \bigsqcup_{a \in A} N_a
\]
of $(N, \pi)$ into complete Poisson submanifolds
and a family of weak sections 
$i_a \colon N_a \hra M$ of $p$ (one for each stratum
of the partition).
\ede

In the category of algebraic varieties we require $M$ and $N$ to be
smooth algebraic varieties, $X$ and $N_a$ to be locally closed smooth algebraic 
subsets, and $p, i, i_a$ to be algebraic maps.

We have:

\bpr{symlv} Consider a surjective
Poisson submersion $p \colon (M, \Pi) \to (N, \pi)$.
Let $N = \bigsqcup_{a \in A} N_a$ and $i_a \colon N_a \colon M$,
$a \in A$ define a weak splitting of $p$. Then 
for all $a \in A$ the following hold:

(i) Every symplectic leaf of $(N_a, \pi|_{N_a})$ has the form 
$i_a^{-1}((N_a \cap S)^\ci)$ where $S$ is a symplectic leaf of
$(M, \Pi)$. 

(ii) Each symplectic leaf $S'$ of $(N_a, \pi|_{N_a})$ is explicitly realized 
as a symplectic submanifold 
\[
i_a \colon (S', \pi_{S'}) \hra (S, \Pi_S)
\]
of the unique symplectic leaf $S$ of $(M, \Pi)$ that contains
$i_a(S')$.
\epr

All weak sections and weak splittings that we construct in this paper
will have the property that their images are Poisson--Dirac submanifolds 
which admit Dirac projections, i.e. the images will satisfy the conditions 
in \prref{PD2}. 
\sectionnew{Weak sections of quotients of Drinfeld doubles}
\label{sect3}
\noindent
We return to the setting of the introduction: Start with a Manin triple $(\d, \g_+, \g_-)$ where 
$\d$ is a quadratic Lie algebra with (a fixed) nondegenerate 
invariant symmetric bilinear form $\lcor.,.\rcor$ and $\g_\pm$
are two Lagrangian subalgebras such that $\d = \g_+ \oplus \g_-$
as vector spaces. Let $D$ be a connected Lie group with Lie algebra $\d$ and 
$G_\pm$ be the connected subgroups of $D$ with Lie algebras $\g_\pm$. 
Fix a pair of dual bases $\{x_j \}$ and $\{\xi_j\}$ 
of $\g_+$ and $\g_-$ with respect to $\lcor . , . \rcor$. 
The standard $r$-matrix $r = \frac{1}{2} \sum \xi_j \wedge x_j$
gives rise to the Poisson structures
\begin{equation}
\label{DH}
\pi = L(r) - R(r) \quad \mbox{and} \quad 
\pi' = L(r) + R(r)
\end{equation}
on $D$. 
Then $(D, \pi)$ is a Poisson--Lie group and 
$G_\pm$ are Poisson--Lie subgroups. Moreover, 
$(D, \pi)$ is a Drinfeld double of $(G_\pm, \pi|_{G_\pm})$
and $(D, \pi')$ is a Heisenberg double of $(G_\pm, \pi|_{G_\pm})$. 
Finally, $(G_-, - \pi|_{G_-})$ is a dual
Poisson--Lie group of $(G_+, \pi|_{G_+})$.

Set for brevity
\begin{equation}
\label{Npm}
N_\pm := N(\g_\pm), \; \n_\pm = \Lie (N_\pm) = \n(\g_\pm).
\end{equation} 
Denote the canonical projections $p_\pm \colon \d \to \g_\pm$
along $\g_\mp$.
Identify $\d^*$ with $\d$ using the form $\lcor.,.\rcor$
and denote by $\al(y)$ the right invariant 1-form 
on $D$ corresponding to $y \in \d \cong \d^*$.  

The bundle maps $\pi^\sharp \colon T^*D \to T D$ and
$(\pi')^\sharp \colon T^* D \to TD$ are given by 
the following formulas.

\ble{sharp} In the above setting, for all 
$x \in \g_+$, $\xi \in \g_-$ and $d \in D$:
\begin{multline}
\label{pi-sh}
\pi^\sharp(\al_d(x+\xi)) = R_d(x) - L_d(p_+ \Ad_d^{-1}(x+\xi) ) \\ =
-R_d(\xi) + L_d(p_- \Ad_d^{-1}(x+\xi) )
\end{multline}
and
\begin{multline}
\label{pi'-sh}
(\pi')^\sharp(\al_d(x+\xi)) = R_d(x) - L_d(p_- \Ad_d^{-1}(x+\xi) ) \\ = 
-R_d(\xi) + L_d(p_+ \Ad_d^{-1}(x+\xi) ).
\end{multline}
\ele
Eq. \eqref{pi'-sh} is proved in \cite{EL2}, eqs. (6.2)--(6.4). 
Eq. \eqref{pi-sh} is analogous.
 
The Poisson structure 
$\pi$ vanishes on 
\begin{equation}
\label{H}
H := N_+ \cap N_- = N(\g_+) \cap N(\g_-),
\end{equation}
see e.g. \cite[Lemma 1.12]{LY}. 
Denote $\h = \Lie \, H$. 
Recall that the right and left actions 
of $(D, \pi)$ on itself and the Heisenberg double
$(D, \pi')$ are Poisson. Thus the left and right action of $H$ 
on $D$ preserves both $\pi$ and $\pi'$. Denote 
their reductions with respect to the right action of $H$
by 
\[
\pi_H \; \; \mbox{and} \; \; \pi'_H \in \Ga( D/H, \wedge^2 T(D/H) ),
\]    
respectively. Thus the canonical projections
\[
\nu \colon (D, \pi) \to (D/H, \pi_H) \; \; 
\mbox{and} \; \; 
\nu' \colon (D, \pi') \to (D/H, \pi'_H)
\]
are Poisson.

Since $\Lie (N_+) \supseteq \g_+$, it follows from the definition of the Drinfeld and Heisenberg 
double Poisson structures \eqref{DH} that
\[
\pi_{D/N_+} := \chi(r)
\]
is a Poisson structure on $D/N_+$ and that the standard projections 
\[
\mu \colon (D, \pi) \to (D/N_+, \pi_{D/N_+})
\quad
\mbox{and} 
\quad
\mu' \colon (D, \pi') \to (D/N_+, \pi_{D/N_+})
\]
are Poisson. 
Denote by 
\[
\eta \colon (D/H, \pi_H) \to (D/N_+, \pi_{D/N_+}) \quad
\mbox{and} \quad
\eta' \colon (D/H, \pi'_H) \to (D/N_+, \pi_{D/N_+})
\]
the induced surjective Poisson submersions. (They are
both Poisson since $\mu = \eta \nu$, $\mu'= \eta' \nu'$,
$\nu$, $\nu'$ are Poisson and 
$\nu$, $\nu'$ are surjective.)

The following proposition will be used in  
our general construction of weak sections for $\eta$:

\bpr{1} Assume that for a given $d \in D$ such that
$dHd^{-1} \subset N_-$ there exists a subgroup $Q$ of 
$D$ with Lie algebra $\q$ satisfying
\begin{align}
&\n_- = \n_- \cap \Ad_d(\n_+) + \n_- \cap \Ad_d(\q) \; \; \mbox{and}
\label{int1}
\\
&Q \cap N_+ = H, \; \n_+ + \q= \n_+ + \q^\perp + \Ad_d^{-1}(\n_-) = \g.  
\label{int2}
\end{align}
Set 
\begin{equation}
\label{Gd}
G_d = N_- \cap d Q d^{-1}.
\end{equation}
Then 
\begin{equation}
\label{wtE}
\wt{E}^d \to G_d d, \; \; 
\wt{E}_{gd}^d := R_{gd}(p_+ \Ad_d(\q^\perp)) + L_{gd}(\n_+), \; \; 
\forall g \in G_d
\end{equation} 
is a subbundle of $T_{G_d d} D$ such that
\begin{equation}
\label{Econd}
\wt{E}^d \cap T_{G_d d}= L(\h), \; \; \wt{E}^d + T (G_d d) = T_{G_d d} D,
\end{equation}
and 
\[
\pi_{G_d d} \in \Ga(G_d d , \wedge^2 \wt{E}^d + \wedge^2 T(G_d d)).  
\]
\epr
\noindent
In \eqref{Econd} $L(\h)$ denotes the subbundle of $T_{G_d d}$
spanned by left invariant vector fields $L(h)$, $h \in \h$.
Here and below $(.)^\perp$ refers to the orthogonal complement 
with respect to $\lcor.,.\rcor$.

The condition $dHd^{-1} \subset N_-$ ensures that $G_d d$ is 
stable under the right action of $H$. The subbundle $\wt{E}^d$ of 
$T_{G_d d} D$ is equivariant with respect to this action. 
Indeed, if $h' \in H$, then 
\begin{multline}
R_{h'} \wt{E}_{gd}^d = R_{gdh'}( p_+ \Ad_d(\q^\perp)) 
+ L_{gd} R_{h'} (\n_+) \\
= R_{gdh'}( p_+ \Ad_d(\q^\perp)) + L_{gdh'} (\n_+) = \wt{E}_{gdh'}^d,
\end{multline}
because $h' \in H \subset N_-$. Therefore, 
the pushforward of $\wt{E}^d$ to $G_d \cdot d H$ is a 
subbundle of $T_{G_d \cdot d H} (D/H)$. As an immediate 
consequence of \prref{1} we obtain the following Corollary:

\bco{1} If, in the above setting, a subgroup $Q$ of $D$ and $d \in D$
satisfy \eqref{int1}-\eqref{int2} and $d H d^{-1} \subset N_-$,
then the submanifold $G_d \cdot d H$ of the quotient $(D/H, \pi_H)$ 
of the Drinfeld double $(D, \pi)$ is a Poisson--Dirac submanifold admitting 
a Dirac projection with associated vector bundle equal to the pushforward 
$E^d= \nu_*(\wt{E}^d)$ of $\wt{E}^d$ to $G_d \cdot dH$.
\eco
\noindent
{\em{Proof of \prref{1}.}} Throughout the proof $g$ 
will denote an element of $G_d$.

First we prove that $\wt{E}^d$ is a subbundle of 
$T_{G_d d} D$ and that \eqref{Econd} holds.
Fix $g \in G_d$. We have:
\begin{align*}
T_{gd} (G_d d) + \wt{E}_{gd}^d  
&\supset L_{gd}( \Ad_d^{-1} \Ad_g^{-1}(\n_- \cap \Ad_d(\q)) + 
L_{gd}(\n_+) 
\\
&= L_{gd}(\Ad_d^{-1}(\n_-)\cap \q + \n_+) 
\supset L_{gd} (\Ad^{-1}_{d}(\n_-)) 
\\
&= L_{gd} (\Ad^{-1}_{d} \Ad_{g}^{-1}(\n_-)) = R_{gd}(\n_-).
\end{align*}
The second inclusion  in the chain follows from \eqref{int1}. Thus,
\begin{align*}
&T_{gd} (G_d d) + \wt{E}_{gd}^d  
 \supset R_{gd}(\n_- + p_+(\Ad_d(\q^\perp)) + L_{gd}(\n_+) \\
&\supset R_{gd}( \n_- + \Ad_d (\q^\perp) ) + L_{gd}(\n_+) \\
&= L_{gd}( \Ad_d^{-1} \Ad_g^{-1} (\n_-) + \Ad_d^{-1} \Ad_g^{-1} \Ad_d (\q^\perp) ) 
+ L_{gd}(\n_+)
\\
&= L_{gd}(\n_+ + \q^\perp + \Ad_d^{-1}(\n_-)) = T_{gd}D,
\end{align*}
where we used \eqref{int1}--\eqref{int2} and the fact 
that $Q$ normalizes $\q^\perp$. Clearly, 
\[
T_{gd} (G_d d) \supset L_{gd} (\h) \; \; \; \mbox{and} \; \; \;  
\wt{E}_{gd}^d \supset L_{gd} (\h).
\]
We claim that
\begin{equation}
\label{dims}
\dim(p_+ \Ad_d(\q^\perp)) + \dim \n_+ + 
\dim (\n_- \cap \Ad_d(\q)) = \dim \g + \dim \h.
\end{equation}
This implies that
$\wt{E}^d$ is a subbundle of $T_{G_d d} D$ and 
that the first equality of \eqref{Econd} is satisfied.
It also follows from \eqref{dims} that $\wt{E}^d$ is the 
direct sum of the subbundles $R(p_+ \Ad_d(\q^\perp))$
and $L(\n_+)$ of $T_{G_d d} D$.

Since $\n_- = \g_- + \h$ and 
$\Ad_d^{-1}(\h) \subset \h \subset \q$,
we have
\begin{align*}
&\dim (\n_- \cap \Ad_d(\q)) = 
\dim \n_- - \dim \g_- + \dim (\g_- \cap \Ad_d(\q))
\\
&= \dim \n_- -\dim \g_- + \dim \g -
\dim( \g_-^\perp + \Ad_d(\q^\perp) )
\\
&=\dim \n_- - \dim(p_+\Ad_d(\q^\perp)).
\end{align*}
Taking into account that 
$\dim \n_+ + \dim \n_- = \dim \d + \dim \h$
leads to \eqref{dims}.

Since 
\[
T_{gd} G_d d = R_{gd} (\n_- \cap \Ad_d(\q))
\]
and 
\[
(\n_- \cap \Ad_d(\q))^\perp = \n_-^\perp + \Ad_d(\q^\perp)
\subset \g_- + p_+(\Ad_d(\q^\perp)),
\]
we have 
\[
(T_{gd} G_d d)^0 \subset \{\al_{gd}(x+\xi) \mid x \in p_+(\Ad_d(\q^\perp)), 
\xi \in \g_- \}.
\]   
Applying \eqref{pi-sh} gives
\[
\pi^\sharp ((T_{gd} G_d d)^0) \subset 
R_{gd}(p_+(\Ad_d(\q^\perp)) + L_{gd}(\g_+)\subset \wt{E}_{gd}^d.
\]
\qed

Observe that for $d \in D$ the conditions \eqref{int1}--\eqref{int2}
ensure that the product $(N_- \cap dQ d^{-1})(N_- \cap dN_+d^{-1})$
is open in $N_-$. This implies that 
$(N_- \cap dQ d^{-1}) \cdot dQ$ is open in $N_- \cdot dQ$
which is a complete Poisson submanifold of $(D/Q, \pi_{D/Q})$,
\cite[Theorem 2.3]{LY}. Thus, \eqref{int1}--\eqref{int2} 
imply that $(N_- \cap dQ d^{-1}) \cdot dQ$ 
is a Poisson submanifold of $(D/Q, \pi_{D/Q})$.

\bth{main1} Assume that for a given $d \in D$ 
such that $d H d^{-1} \subset N_-$ there exists
a subgroup $Q$ of $D$ satisfying the conditions 
\eqref{int1}-\eqref{int2}. Then the smooth map
$i \colon G_d \cdot d N_+ \to D/H$ defined
by $i(g d N_+) = gdH$ for $g \in G_d := N_- \cap d Qd^{-1}$ 
is a weak section of the surjective Poisson submersion
$\eta \colon (D/H, \pi_H) \to (D/N_+, \pi_{D/N_+})$
over $G_d \cdot d N_+$. Its image is 
a Poisson--Dirac submanifold of $(D/H, \pi_H)$
admitting a Dirac projection with associated 
bundle $E^d := \nu_*(\wt{E}^d)$ where $\wt{E}^d$ is given by
\eqref{wtE}.
\eth
\begin{proof}
It is straightforward to check that $i$ is well defined:
If $g_1, g_2 \in G_d$ and $g_1 d Q = g_2 d Q$, 
then $(g_2)^{-1} g_1 \in N_- \cap d (Q \cap N_+) d^{-1}
\subset d H d^{-1}$ because of \eqref{int2}. Thus 
$g_1 d H = g_2 d H$.

For $g \in G_d$, \coref{1} implies that 
\[
(\pi_H)_{gdH} \in \wedge^2 T_{gdH} (G_d \cdot d H) + 
\wedge^2 E_{gdH}^d.   
\]
Observe that $E_{gdH}^d$ contains the tangent space 
$\nu_*(L_{gd}(\n_+))$ to the fiber of $\nu$ through
$gdH$. Since $\nu_* (\pi_H) = \pi_{D/N_+}$ 
and $\nu \circ i = \id_{G_d \cdot d H}$, we have that 
the projection of $(\pi_H)_{gdH}$ to 
$\wedge^2 T_{gdH} (G_d \cdot d H)$ is 
$i_*(\pi_{D/N_+}|_{G_d \cdot dH})$.
This completes the proof of Theorem.
\end{proof}

The next theorem provides a sufficient condition for the existence of 
a weak splitting of the surjective Poisson submersion in 
\thref{main1}.

\bth{main2} Assume that $Q$ is a subgroup of $D$ with Lie algebra
$\q$ such that there exists a set of representatives
$\DD$ for the $(N_-,N_+)$-double cosets of $D$ 
satisfying
\begin{align}
&N_- = (N_- \cap dQd^{-1})(N_- \cap dN_+d^{-1})\; \; \mbox{and}
\label{int3}
\\
&dH d^{-1} \subset N_-, \; Q \cap N_+ = H, \; \n_+ + \q= 
\n_+ + \q^\perp + \Ad_d^{-1}(\n_-) = \g
\label{int4}
\end{align}
for all $d \in \DD$.
Then the partition 
\begin{equation}
\label{part}
D/N_+ = \sqcup_{d \in \DD} (N_- \cap dQd^{-1}) \cdot d N_+
\end{equation}
and the family of smooth maps 
\[
i_d \colon (N_- \cap dQd^{-1}) \cdot d N_+ \to D/H, \; \; \;
i_d(g d N_+) := gdH, \; \; 
\forall g \in N_- \cap d Qd^{-1}
\]
is a weak splitting of the surjective Poisson submersion
\[
\eta \colon (D/H, \pi_H) \to (D/N_+, \pi_{D/N_+}).
\]
In addition, the images of $i_d$ are  Poisson--Dirac submanifolds 
of $(D/H, \pi_H)$ admitting Dirac projections with associated 
bundles $E^d := \nu_*(\wt{E}^d)$ for the bundles $\wt{E}^d$ given by
\eqref{wtE}.
\eth

\begin{proof}
The condition \eqref{int3} implies that 
\[
(N_- \cap dQd^{-1}) \cdot d N_+ = N_- \cdot d N_+.
\]
It follows from the definition of the set $\DD$ 
that \eqref{part} defines a partition of $G/N_+$.
This equality also implies that each stratum of the partition 
is a complete Poisson submanifold of $(D/N_+, \pi_{D/N_+})$,
because this is a property of all $N_+$-orbits on $D/N_+$.
The rest of the theorem 
follows from \thref{main1}.
\end{proof}

\prref{symlv} (i) implies that in the setting of \thref{main1} 
each symplectic leaf of $(D/N_+, \pi_{D/N_+})$ inside the stratum 
$(N_- \cap dQd^{-1}) \cdot d N_+$ is of the form 
\[
i_d^{-1}(( i_d(N_- \cap dQd^{-1}) \cap S)^\circ) 
\]
for a symplectic leaf of $S$ of $(D/H, \pi_H)$. Furthermore, 
by \prref{symlv} (ii) each symplectic leaf of 
$(D/N_+, \pi_{D/N_+})$ is explicitly realized 
as a symplectic submanifold of a symplectic leaf of 
$(D/H, \pi_H)$ via one of the maps $i_d$.

\bre{conn} In \prref{1}, and Theorems \ref{tmain1} and \ref{tmain2} one can replace
$N_\pm$ with any pair of subgroups $N'_\pm$ of $D$ such that 
$N_\pm^\circ \subset N'_\pm \subset N_\pm$. The corresponding 
statements hold true for $H := N_+' \cap N'_-$.
Their proofs are analogous and are left to the reader.
\ere

\bre{algD} If the group $D$ and the subgroup $Q$ in Theorems \ref{tmain1} 
and \ref{tmain2} are algebraic, then the constructed weak sections and splittings 
are algebraic.
\ere
\sectionnew{Weak sections of quotients of Heisenberg doubles}
\label{sect4}
\noindent
In this section we prove results for quotients of Heisenberg doubles 
that are similar to the results from the previous section for 
quotients of Drinfeld doubles. We use the setting and notation 
of the previous section.
Using the second part of \leref{sharp} one proves the 
following analog of \coref{1} and \thref{main1} for Heisenberg
doubles. We omit its proof since it is analogous to the case of Drinfeld 
doubles.

\bth{main1H} Let $d$ and a subgroup $Q$ of $D$ 
with Lie algebra $\q$ satisfy
\begin{align}
&\n_+ = \n_+ \cap \Ad_d(\n_+) + \n_+ \cap \Ad_d(\q) \; \; \mbox{and}
\label{int1H}
\\
&dHd^{-1} \subset N_+, \; Q \cap N_+ = H, \; \n_+ + \q = 
\n_+ + \q^\perp + \Ad_d^{-1}(\n_+) = \g.
\label{int2H}
\end{align}
Set $G'_d = N_+ \cap d Q d^{-1}$. Then the submanifold $G'_d \cdot d H$ 
of the quotient $(D/H, \pi'_H)$ of the Heisenberg double $(D, \pi')$ is a 
Poisson--Dirac submanifold admitting a Dirac projection with associated 
vector bundle equal to the the pushforward $F^d:= \nu'_*(\wt{F}^d)$ 
of the vector bundle
\begin{equation}
\label{Fbundl}
\wt{F} \to G'_d d, \; \;
\wt{F}_{gd} := R_{gd}(p_- \Ad_d(\q^\perp)) + L_{gd}(\n_+), \; \;
\forall g \in G'_d.
\end{equation}
In addition, the map $i \colon G'_d \cdot d N_+ \to D/H$ defined
by $i(g d N_+) = gd H$ for $g \in G'_d$
is a weak section of the surjective Poisson submersion
$\eta \colon (D/H, \pi'_H) \to (D/N_+, \pi_{D/N_+})$
over $G'_d \cdot d N_+$.
\eth

As in the previous section the theorem implies the following:

\bth{main2H} Let $Q$ be a subgroup of $D$ with Lie algebra
$\q$ for which there exists a set of representatives
$\DD \subset N(H)$ for the $(N_+,N_+)$-double cosets of $D$
satisfying
\begin{equation}
N_+ = (N_+ \cap dQd^{-1})(N_+ \cap dN_+d^{-1})
\label{int3H}
\end{equation}
and \eqref{int2H} for all $d \in \DD$. Then the partition
\[
D/N_+ = \sqcup_{d \in \DD} (N_+ \cap dQd^{-1}) \cdot d N_+
\]
and the family of maps
\[
i'_d \colon (N_+ \cap dQd^{-1}) \cdot d N_+ \to D/H, \; \; \;
i'_d(g d N_+) := gdH, \; \; \forall
g \in N_+ \cap d Qd^{-1}
\]
provide a weak splitting of the surjective Poisson submersion
$\eta' \colon (D/H, \pi'_H) \to (D/N_+, \pi_{D/N_+})$.
In addition, the images of $i'_d$ are  Poisson--Dirac submanifolds 
of $(D/H, \pi'_H)$ admitting Dirac projections with associated 
bundles $F^d := \nu'_*(\wt{F}^d)$ for the bundles $\wt{F}^d$ given by
\eqref{Fbundl} with $G'_d = N_+ \cap dQd^{-1}$.
\eth

In light of \prref{symlv}, \thref{main2H} provides an explicit 
realization of the symplectic leaves of $(D/N_+, \pi_{D/N_+})$ 
as symplectic submanifolds of the symplectic leaves of
the Heisenberg double $(D/H, \pi'_H)$.

As in the case of Drinfeld doubles, in Theorems \ref{tmain1H} and \ref{tmain2H} 
one can replace $N_\pm$ with any pair of subgroup $N'_\pm$ of $D$ such that 
$N_\pm^\circ \subset N'_\pm \subset N_\pm$ in which case one sets 
$H = N'_+ \cap N'_-$. The proofs of those slightly more general 
statements are analogous.

If the groups $D$ and $Q$ are algebraic, the above constructed weak sections 
and weak splittings are also algebraic.
\sectionnew{Applications to flag varieties}
\label{sect5}
\noindent
This section contains applications of the results from the previous two sections to 
Poisson structures on flag varieties. Subsections \ref{5.1}--\ref{5.3} deal with 
complex algebraic Poisson structures. There we construct weak splittings 
for complex surjective Poisson submersions from Drinfeld and Heisenberg doubles
to flag varieties, double flag varieties and certain natural multi-flag generalizations. 
In \S \ref{5.4} we give applications to real algebraic Poisson structures
on flag varieties. 

We note that all of the weak splittings that are constructed in this section 
provide (via \prref{symlv}) explicit realizations of the symplectic leaves 
of Poisson structures on flag varieties 
as symplectic submanifolds of symplectic leaves of Drinfeld and Heisenberg doubles.

Throughout the section $G$ will denote an arbitrary connected complex simple Lie group
with Lie algebra $\g$. The Killing form on $\g$ will be denoted 
by $\lcor.,.\rcor$. We fix a pair of opposite 
Borel subgroups $B_\pm$ of $G$ and 
the corresponding maximal torus $T := B_+ \cap B_-$. 
Let $U_\pm$ be the unipotent radicals
of $B_\pm$. Set 
\[
\h := \Lie \, T, \; \; 
\b_\pm := \Lie B_\pm, \; \; 
\mbox{and} \; \;  
\uu_\pm := \Lie U_\pm.
\]

Denote the Weyl group of $G$ by $W$. For 
each $w \in W$, fix a representative $\dot{w}$ in the 
normalizer of $T$ in $G$.
\noindent
\subsection{}
\label{5.1}
The simplest applications of Theorems \ref{tmain1} and \ref{tmain1H} 
are to weak splittings for Poisson structures on 
flag varieties.

Recall the standard Manin triple
\[
\d := \g \oplus \h, \;
\g_\pm := \{ (x_\pm +h, \pm h) \mid x_\pm \in \uu_\pm, h \in \h \}
\]
with respect to the invariant bilinear form on $\d$ given by
\[
\lcor (x_1,y_1), (x_2,y_2) \rcor := 
\lcor x_1, x_2 \rcor - \lcor y_1, y_2 \rcor, \; \; 
\forall x_i \in \g, y_i \in \h.
\]
The Drinfeld and Heisenberg double Poisson structures on $D:=G \times T$ 
will be denoted by $\pi$ and $\pi'$, respectively. In this case 
$N_\pm = N(\g_\pm) = B_\pm \times T$ 
and $H = N_+ \cap N_- = T \times T$.
The reductions of $\pi$ and $\pi'$ to 
$G/T \cong (G \times T) /(T \times T) = D/H$ 
will be denoted by $\pi_T$ and $\pi'_T$.
Both structures reduce to the same Poisson 
structure on the flag variety 
$G/B_+ \cong (G\times T) /(B_+ \times T) = D / N_+ $ 
called the standard Poisson structure. The latter will be denoted
by $\pi_{G/B_+}$. 

It is easy to verify that the group $Q= N_- = B_- \times T$
and the set $\DD = \{ \dot{w} \mid w \in W\}$ 
satisfy the conditions in \thref{main2} for the above choice of $D$ and $\g_\pm$.
This implies the following result of Goodearl and the author 
proved in \cite[Theorem 3.2]{GY}. 
For its statement we need to introduce some additional notation.
Denote the vector bundle
\[
\wt{E}_w \to \big( (B_- \cap w B_+ w^{-1}) \times T \big) (\dot{w}, 1)
\]
with fibers
\[
\wt{E}^w_g:= R_g \big(p_+ \Ad_{(\dot{w},1)}(\uu_- \oplus 0) \big) 
+ L_g(\b_+ \oplus \h)
\]
for $g \in \big( (B_- \cap w B_+ w^{-1}) \times T \big) (\dot{w}, 1)$
where $p_+ \colon \d \to \g_+$ is the projection along $\g_-$.
Here the direct sum notation is used to denote subspaces of $\g \oplus \h$
identified with $\Lie (G \times T)$. 
Denote the canonical projection
\begin{equation}
\label{nu2}
\nu \colon G \times T \to (G \times T) / (B_+ \times T) \cong G/B_+.
\end{equation}
By \coref{1} the pushforward $E^w:= \nu_*(\wt{E}^w)$   
is a well defined vector bundle over the Schubert cell 
$B_- \cdot w B_+ \subset G/B_+$.

\bth{FlagD} \cite{GY} For all complex simple Lie groups $G$, 
the partition of the full flag variety $G/B_+$ into Schubert cells
\[
G/B_+ = \sqcup_{w \in W} B_- \cdot w B_+
\]
and the family of maps
\[
i_w \colon B_- \cdot w B_+ \to G/T, \; \; 
i_w(b_- \dot{w} B_+) := b_- \dot{w} T, \; \; \forall b_- \in B_- \cap w B_- w^{-1}
\]
define a weak splitting of the surjective Poisson submersion
\[
(G/T, \pi_T) \to (G/B_+, \pi_{G/B_+})
\]
from a Drinfeld double to the flag variety.
Furthermore, the image of each map $i_w$ is a Poisson--Dirac 
submanifold of $(G/T, \pi_T)$ admitting a Dirac projection 
with associated bundle $E^w$ defined above.
\eth

\thref{main2H} for weak splittings of quotients of Heisenberg doubles is also applicable to 
flag varieties to obtain a weak splitting for the surjective 
Poisson submersion 
\[
(G/T, \pi'_T) \to (G/B_+, \pi_{G/B_+}).
\]
It is easy to verify that the conditions 
of \thref{main2H} are satisfied by for the same  
group $Q= N_- = B_- \times T$, set $\DD = \{ \dot{w} \mid w \in W \}$ 
and the current choice of $D$ and $\g_\pm$. 
Applying the theorem leads to the following result:
  
\bth{FlagH} For all complex simple Lie groups $G$, the partition of the full flag 
variety $G/B_+$ into Schubert cells
\[
G/B_+ = \sqcup_{w \in W} B_+ \cdot w B_+ 
\]
and the family of maps
\[
i'_w \colon B_+ \cdot w B_+ \to G/B_+, \; \;
i'_w(b_+ w B_+) = b_+ \dot{w} T, \; \; b_+ \in B_+ \cap w B_- w^{-1}
\]
is a weak splitting of the surjective Poisson submersion
\[
(G/T, \pi'_T) \to (G/B_+, \pi_{G/B_+})
\]
from a Heisenberg double to the flag variety.
The image of each map $i'_w$ is a Poisson--Dirac submanifold
of $(G/T, \pi'_T)$ admitting a Dirac projection 
with associated bundle $\nu_*(\wt{F}^w)$ for the pushforward 
bundle $\nu_*(\wt{F}^w)$ with respect to \eqref{nu2} where
\[
\wt{F}^w \to \big( (B_+ \cap w B_- w^{-1}) \times T \big) (\dot{w}, 1)
\]
is the vector bundle with fibers
\[
\wt{F}^w_g:= R_g \big(p_- \Ad_{(\dot{w},1)}(\uu_- \oplus 0) \big) 
+ L_g(\b_+ \oplus \h)
\]
for $g \in \big( (B_+ \cap w B_- w^{-1}) \times T \big) (\dot{w}, 1)$
and $p_- \colon \d \to \g_-$ is the projection along $\g_+$.
\eth

Because of of \prref{symlv}, Theorems \ref{tFlagD} and \ref{tFlagH} provide 
an explicit realization of the symplectic leaves of the flag varieties 
$(G/B_+, \pi_{G/B_+})$ as symplectic submanifolds of the symplectic leaves of
Drinfeld and Heisenberg doubles. 

We note that the partitions into Schubert cells in Theorems \ref{tFlagD} and \ref{tFlagH} are with respect 
to opposite Borel subgroups. The two results can be derived from each other.
Let $w_\circ$ be the longest element of $W$. The equivalence is shown using 
the facts that the translation action of $\dot{w}_\circ$ on $(G/B_+, \pi_{G/B_+})$ is anti-Poisson,
and the left translation action of $\dot{w}_\circ$ on $G \times T$ interchanges $\pi$ and $-\pi'$.
\subsection{}
\label{5.2}
Next, we consider the standard Manin triple
\begin{multline}
\label{Man}
\d := \g \oplus \g, \; \; \g_+ := \{ (x_+ +h, x_- - h) \mid x_\pm \in \uu_\pm, 
h \in \h \}, \\
\g_- := \{(x, x) \mid x \in \g \}
\end{multline}
with respect to the invariant bilinear form on $\d$
\[
\lcor (x_1,y_1), (x_2,y_2) \rcor := 
\lcor x_1, x_2 \rcor - \lcor y_1, y_2 \rcor, \; \; 
x_i, y_i \in \g.
\]
Let $D := G \times G$.
For this setting we have 
\[
N_+ = N(\g_+) = B_+ \times B_-, \; \;
N_-^\circ = N(\g_-)^\circ = G_\Delta, \; \; 
\mbox{and} \; \;  N_+^\circ \cap N_- = T_\Delta 
\]
where $G_\Delta$ and $T_\Delta$ 
are the diagonal subgroups of $G \times G$ and $T \times T$,
respectively.  

The Drinfeld and Heisenberg double Poisson structures on 
$G \times G$ will be again denoted by $\pi$ and $\pi'$.
The group $T_\Delta$ acts on the left and the right 
on $(G \times G, \pi)$ and $(G \times G, \pi')$ by 
Poisson automorphisms.
The reductions of $\pi$ and $\pi'$ to $(G \times G)/T_\Delta$ 
will be denoted by $\pi_{T_\Delta}$ and $\pi'_{T_\Delta}$.
The pushforwards of $\pi_{T_\Delta}$ and $\pi'_{T_\Delta}$
to $(G \times G)/ N_- \cong G/B_+ \times G/B_-$ are well defined 
and are equal to each other. 
Denote the resulting Poisson structure by $\pi_{df}$. 
The Poisson manifold 
\[
(G/B_+ \times G/B_-, \pi_{df})
\]
is the double flag variety from \cite{WY} (up to a rescaling of the Poisson structure).

\thref{main1} cannot be applied to the submersion
$((G\times G) /T_\Delta, \pi_{T_\Delta}) \to
(G/B_+ \times G/B_-, \pi_{df})$, but \thref{main2H} 
can be applied for the following choice of a group $Q$
and a set $\DD$:
\[
Q = (U_- \times U_+) T_\Delta, \quad
\DD = \{ (\dot{w}, \dot{v}) \mid w, v \in W \},
\]
and the above $D$ and $\g_\pm$. 
This gives a weak splitting of the surjective Poison submersion
$((G\times G) /T_\Delta, \pi'_{T_\Delta}) \to
(G/B_+ \times G/B_-, \pi_{df})$. We will need the 
following notation to state the result. 
Denote the projection 
\begin{equation}
\label{nu-d}
\nu' \colon G \times G \to (G \times G) / T_\Delta.
\end{equation}
Consider the vector bundle 
\[
\wt{F}^{w,v} \to \big( (U_+ \cap w U_+ w^{-1}) \times (U_- \cap v U_- v^{-1}) \big) \times T_\Delta
\]
with fibers
\begin{equation}
\label{buF}
\wt{F}^{w,v}_g:= R_g(p_-\Ad_{\dot{w}, \dot{v}}
( \uu_- \oplus \uu_+ + \tt_{a \De} ) + L_g(\b_+ \oplus \b_-)
\end{equation}
for $g \in \big( (U_+ \cap w U_+ w^{-1}) \times (U_- \cap v U_- v^{-1}) \big) \times T_\Delta$
where the direct sum notation is used for subspaces of $\d= \g \oplus \g$ identified with 
$\Lie (G \times G)$, $p_- \colon \d \to \g_-$ denotes the projection along $\g_+$,
and $\tt_{a \De}$ is the antidiagonal of $\tt \oplus \tt$.

\bth{FlagH2} For all complex simple Lie groups $G$, the partition
of the double flag variety into $B_+ \times B_-$-Schubert cells
\[
G/B_+ \times G/B_- = \sqcup_{w,v \in W} 
B_+ \cdot w B_+ \times B_- \cdot v B_-
\]
and the family of maps
\[
i'_{w,v} \colon B_+ \cdot w B_+ \times B_- \cdot v B_-
\to (G \times G)/T_\Delta
\]
given by
\[
i'_{w,v} (u_+ w B_+, u_- v B_-) := (u_- \dot{w}, u_+ \dot{v})  T_\Delta
\]
for all $u_+ \in U_+ \cap w U_+ w^{-1}$ and 
$u_- \in U_- \cap v U_- v^{-1}$
is a weak splitting of the surjective Poisson submersion
\[
((G \times G) / T_\Delta , \pi'_\Delta) \to (G/B_+ \times G/B_-, \pi_{df})
\]
from a Heisenberg double to the double flag variety. The image 
of each map $i'_{w,v}$ is a Poisson--Dirac submanifold of 
$((G \times G) / T_\Delta , \pi'_\Delta)$ admitting a Dirac projection 
with associated vector bundle $\nu'_*(\wt{F}^{w,v})$, cf.
\eqref{nu-d} and \eqref{buF}.
\eth
\subsection{}
\label{5.2b} All results in Sect. \ref{sect5} on weak splittings for 
full flag varieties have analogs to 
partial flag varieties. 
We will provide full details in the case 
of double partial flag varieties. The generalizations of the results in 
\S \ref{5.1}, \ref{5.3}, and \ref{5.4} are analogous.

For a subset $I$ of simple roots of $\g$ denote by $P^I_\pm \supseteq B_\pm$ the 
corresponding parabolic subgroups of $G$ and by $W_I$ the subgroup of 
the Weyl group. Let $W^I$ be the sets of
minimal length representatives for the cosets in $W/W_I$.  

Fix two subsets $I_1, I_2$ of simple roots of $\g$. The pushforward 
of $\pi_{df}$ under the canonical projection 
\[
G/B_+ \times G/B_- \to G/P^{I_1}_+ \times G/P^{I_2}_- 
\]
is a well defined Poisson structure since $P^{I_1}_+ \times P^{I_2}_-$ 
is a Poisson--Lie subgroup of $(G \times G, \pi)$. Denote this 
pushforward by $\pi^{I_1, I_2}_{df}$. The map 
\[
(G/B_+ \times G/B_-, \pi_{df}) \to (G/P^{I_1}_+ \times G/P^{I_2}_-, \pi^{I_1, I_2}_{df}) 
\]
is a surjective Poisson submersion and its restrictions 
\[
(B_+ \cdot w B_+ \times B_- \cdot v B_-, \pi_{df}) \to 
(B_+ \cdot w P^{I_1}_+ \times B_- \cdot v P^{I_2}_-, \pi^{I_1, I_2}_{df})
\]
are Poisson isomorphisms for all $w \in W^{I_1}$, $v \in W^{I_2}$. 
(Similar Poisson isomorphisms are constructed in the settings of \S \ref{5.1}, \ref{5.3}, and \ref{5.4}.
This produces the generalizations of those results to the cases of partial 
flag varieties.)
Taking inverses of the above Poisson isomorphisms and 
composing them with the maps $i'_{w,v}$ in \thref{FlagH2} 
leads to the following result:

\bco{pFlagH2} For all connected complex simple Lie groups $G$ and subsets of simple roots 
$I_1$ and $I_2$, the partition of the corresponding 
double partial flag variety into $B_+ \times B_-$-Schubert cells
\[
G/P^{I_1}_+ \times G/P^{I_2}_- = \sqcup_{w \in W^{I_1}, v \in W^{I_2}} 
B_+ \cdot w P^{I_1}_+ \times B_- \cdot v P^{I_2}_-
\]
and the family of maps
\[
j'_{w,v} \colon B_+ \cdot w P^{I_1}_+ \times B_- \cdot v P^{I_2}_-
\to (G \times G)/T_\Delta
\]
given by
\[
j'_{w,v} (u_+ w P^{I_1}_+, u_- v P^{I_2}_-) := (u_- \dot{w}, u_+ \dot{v})  T_\Delta
\]
for all $u_+ \in U_+ \cap w U_+ w^{-1}$ and 
$u_- \in U_- \cap v U_- v^{-1}$
is a weak splitting of the surjective Poisson submersion
\[
((G \times G) / T_\Delta , \pi'_\Delta) \to (G/P^{I_1}_+ \times G/P^{I_2}_-, \pi_{df}^{I_1,I_2})
\]
from a Heisenberg double to the double partial flag variety. The image 
of each map $j'_{w,v}$ (which is the same as the image of the map $i'_{w,v}$)
is a Poisson--Dirac submanifold of 
$((G \times G) / T_\Delta , \pi'_\Delta)$ admitting a Dirac projection 
with associated vector bundle $\nu'_*(\wt{F}^{w,v})$, see
\eqref{nu-d} and \eqref{buF}.
\eco
\subsection{}
\label{5.3}
The results in \S \ref{5.1}--\ref{5.2} can be generalized 
to a very large class of Poisson structures on multiple flag varieties.
Those are Cartesian products of flag varieties for complex simple Lie groups 
(i.e., flag varieties for reductive Lie groups) with Poisson structures 
which in general are not products of Poisson structures on the factors. 
Since the arguments are similar, we only state the results leaving 
the details to the reader. 

We start with any reductive Lie algebra $\d$ and an invariant bilinear 
form $\lcor . , . \rcor$ on it. All Lagrangian subalgebras and 
Manin triples in this situation were classified by Delorme \cite{De} 
up to the action of the adjoint group of $\d$.
Let $\ol{\b}_\pm$ are a pair of opposite Borel 
subalgebras of $\d$ and $\ol{\tt} := \ol{\b}_+ \cap \ol{\b}_-$ be the corresponding
Cartan subalgebra of $\d$. Denote by $\ol{\uu}_\pm$ the nilradicals of $\ol{\b}_\pm$.
Let $D$ be a connected reductive Lie group with Lie algebra $\d$, and $\ol{B}_\pm$ and $\ol{T}$ be its 
Borel subgroups and maximal torus corresponding to $\ol{\b}_\pm$ and $\ol{\tt}$.

Consider any Manin triple
\begin{equation}
\label{Man-dd}
(\d, \g_+, \g_-) \; \; \;  \mbox{such that} \; \; \; 
\g_+ \subset \ol{\b}_+.
\end{equation}
The results of Delorme imply that after a conjugation by an element of $B_+$, one has
\[
\g_+ = \uu_+ + \g_+ \cap \ol{\tt}
\] 
and
\begin{equation}
\label{H-dd}
H := N_+ \cap N_-^\circ = N(\g_+) \cap N(\g_-)^\circ = \ol{T} \cap N(\g_-)^\circ.
\end{equation}
In particular, $N_+ = \ol{B}_+$. 
One can also write an explicit formula for the connected 
group $H$ in terms of generalized Belavin--Drinfeld triples in the setting of \cite{De}.
In the rest we will assume that 
the conjugation by an element of $\ol{B}_+$ is performed 
so that the above conditions are satisfied.

Denote the Drinfeld and Heisenberg double Poisson structures on $D$ 
corresponding to a Manin triple of the type \eqref{Man-dd} by 
$\pi$ and $\pi'$. 
By the general facts in Sect. \ref{sect3}, the left and right 
regular actions of $H$ on $(G, \pi')$ are Poisson. Denote the 
reduction $(D/H, \pi'_H)$ for the right action and the surjective 
Poisson submersion
\begin{equation}
\label{nu-dd}
\nu' \colon (D, \pi') \to (D/H, \pi'_H).
\end{equation}
The Drinfeld and Heisenberg
Poisson structures $\pi$ and $\pi'$ descend to the same Poisson structure
on the multiple flag variety $D/N_+ = D/\ol{B}_+$ which will be denoted 
by $\pi_{D/\ol{B}_+}$. The Poisson structures in \S \ref{5.1}--\ref{5.2} are 
special cases of this construction when
\[
D=G \times T \; \; \; \mbox{or} \; \; \; D = G \times G
\]
for a complex simple Lie group $G$ and a maximal torus $T$ of $G$.

The canonical projection
\[
(D/H, \pi'_H) \to (D/\ol{B}_+, \pi_{D/\ol{B}_+})
\]
is Poisson, because $\ol{B}_+ = N(\g_+)$ is a Poisson--Lie subgroup of $(D, \pi)$.
Denote the connected subgroups of $D$ with Lie algebras $\ol{\uu}_\pm$ by $\ol{U}_\pm$. 
Let $\ol{W}$ be the Weyl group of $D$. For each of $w \in \ol{W}$, fix a representative $\dot{w}$  
in the normalizer of the maximal torus $\ol{T}$ of $D$.

A simple computation shows that the conditions of \thref{main2H} are satisfied 
for the group $Q = H \ol{U}_-$ and the set $\DD = \{ \dot{w} \mid w \in \ol{W} \}$.
We have:

\bth{MFlagH} For all Manin triples for a connected reductive algebraic group $D$ of the form 
\eqref{Man-dd}, the partition of the multiple flag 
variety $D/\ol{B}_+$ into Schubert cells
\[
D/\ol{B}_+ = \sqcup_{w \in \ol{W}} \ol{B}_+ \cdot w \ol{B}_+ 
\]
and the family of maps
\[
i'_w \colon \ol{B}_+ \cdot w \ol{B}_+ \to D/H, \; \;
i'_w(u_+ w \ol{B}_+) := u_+ \dot{w} H, \; \; \forall u_+ \in \ol{U}_+ \cap w \ol{U}_- w^{-1}
\]
(recall \eqref{H-dd}) 
is a weak splitting of the surjective Poisson submersion
\[
(D/H, \pi'_H) \to (D/\ol{B}_+, \pi_{D/\ol{B}_+})
\]
from a Heisenberg double to the multiple flag variety.
The image of each map $i'_w$ is a Poisson--Dirac submanifold
of $(D/H, \pi'_H)$ admitting a Dirac projection 
with associated bundle $\nu'_*(\wt{F}^w)$ for the pushforward 
with respect to \eqref{nu-dd} of the vector bundle
\[
\wt{F}^w \to (\ol{U}_+ \cap w \ol{U}_- w^{-1}) \dot{w} H
\]
with fibers
\[
\wt{F}^w_g:= R_g \big(p_- \Ad_{\dot{w}}(\ol{\uu}_- + \h^\perp) \big) 
+ L_g(\ol{\b}_+)
\]
for $g \in (\ol{U}_+ \cap w \ol{U}_- w^{-1}) \dot{w} H$.
Here $\h^\perp$ denotes the orthogonal complement to $\h:= \Lie H$ in $\ol{\tt}$ 
with respect to $\lcor.,.\rcor$,
and $p_- \colon \d \to \g_-$ is the projection along $\g_+$.
\eth
\subsection{}
\label{5.4}
All results in \S \ref{5.1}--\ref{5.3} remain valid when all
complex groups are replaced with their real split forms. These 
provide many examples of weak splittings of surjective 
Poisson submersions to real flag varieties.

We continue with 
certain non-split analogs of the results in \S \ref{5.1} and \ref{5.3}
which concern the real Poisson structures on complex flag varieties introduced 
by Foth and Lu in \cite{FL}. Consider $\g$ as a real quadratic Lie algebra with the 
nondegenerate bilinear form 
\begin{equation}
\label{Im-K}
x, y \in \g \mt \Im \lcor x, y \rcor \in \Rset.
\end{equation}
Each Vogan diagram $v$ for $\g$ gives rise to a complex conjugate linear 
involution $\tau_v$ of $\g$ and to the real form 
$\g_v:= \g^{\tau_v}$ of $\g$, see \cite{K,FL} for details. The map $\tau_v$ 
preserves the complex Cartan subalgebra $\tt$ of $\g$. The following is a (real) Manin triple:
\[
(\d:=\g, \g_+ := \tt^{-\tau_v} + \uu_+, \g_- := \g^v)
\] 
see \cite[Sect. 2]{FL}. Clearly, 
\[
N_+= N(\g_+) = B_+ \quad \mbox{and} 
\quad N_-^\circ = N(\g_-) = G_v 
\]
where $G_v$ is the real form of $G$ associated to $\g_v$. 
The group $H := N_+ \cap N_-^\circ$ is the connected subgroup of $T$ with Lie algebra 
\[
\Lie H = \tt^{-\tau_v}.
\]
Denote once again the associated Drinfeld and Heisenberg double Poisson structures 
on $D:=G$ by $\pi_v$ and $\pi'_v$. We have $\pi|_H=0$, and thus $H$ acts by Poisson automorphisms
on $(G,\pi')$ on the left and right. Denote the reduced Poisson structure 
on $G/H$ by $\pi'_{H,v}$ and the Poisson projection 
\begin{equation}
\label{nu-real}
\nu' \colon (G, \pi'_v) \to (G/H, \pi'_{H,v}).
\end{equation}
The pushforwards of $\pi$ and $\pi'$ under the canonical projection
$G \to G/B_+$ are well defined and are equal to each other, because 
$B_+= N(\g_+)$. Denote the corresponding real Poisson structure 
on the complex flag variety $G/B_+$ by $\pi_{G/B_+,v}$. It has very 
interesting properties, for example the intersections of the orbits 
of the Borel subgroup $B_+$ and the real form $G_v$ are regular complete Poisson 
submanifolds. \thref{main2H} is applicable to construct a 
weak splitting of the real surjective Poisson submersion 
\[
(G/H, \pi'_{H,v}) \to (G/B_+, \pi_{G/B_+,v}).
\] 
Once again it is easy to verify that 
the group $Q= H U_-$ and the set $\DD= \{ \dot{w} \mid w \in W \}$ satisfy the 
conditions of \thref{main2H}. This leads to the following result:

\bth{FlagVogan} For all Vogan diagrams for a complex simple group $G$, 
the partition of the flag variety $G/B_+$ into Schubert cells
\[
G/B_+ = \sqcup_{w \in W} B_+ \cdot w B_+ 
\]
and the family of maps
\[
i'_w \colon B_+ \cdot w B_+ \to G/H, \; \;
i'_w(u_+ w B_+) = u_+ \dot{w} H, \; \; u_+ \in U_+ \cap w U_- w^{-1}
\]
is a weak splitting of the real surjective Poisson submersion
\[
(G/H, \pi'_{H,v}) \to (G/B_+, \pi_{G/B_+,v})
\]
from a Heisenberg double to the flag variety.
The image of each map $i'_w$ is a Poisson--Dirac submanifold
of $(G/H, \pi'_H)$ admitting a Dirac projection 
with associated bundle $\nu'_*(\wt{F}^w)$ for the pushforward 
with respect to \eqref{nu-real} of the vector bundle
\[
\wt{F}^w \to (U_+ \cap w U_- w^{-1}) \dot{w} H
\]
with fibers
\[
\wt{F}^w_g:= R_g \big(p_- \Ad_{\dot{w}}(\uu_- + (\tt^{-\tau_v})^\perp) \big) 
+ L_g(\b_+)
\]
for $g \in (U_+ \cap w U_- w^{-1}) \dot{w} H$.
Here $(\tt^{-\tau_v})^\perp$ denotes the orthogonal complement of $\tt^{-\tau_v}$
in $\tt$ with respect to \eqref{Im-K}
and $p_- \colon \g \to \g_-$ is the projection along $\g_+$.
\eth

\end{document}